\documentclass[runningheads]{llncs}

\usepackage[T1]{fontenc}
\usepackage{lmodern}
\usepackage{setspace}
\usepackage{parskip}
\usepackage{hyperref}

\begin{document}

\title{Ice cream, chocolate, random identification codes, quantum error correction –\\
and a truly warm heart.\\[2mm]
\large In Memory of Vladimir (Volodia) Sidorenko}

\titlerunning{In Memory of Vladimir Sidorenko}

\author{Christian Deppe}
\authorrunning{C. Deppe}

\institute{Technische Universität Braunschweig,\\
Institute for Communications Technology (IfN), Germany\\
\email{christian.deppe@tu-braunschweig.de}}

\maketitle

\onehalfspacing

\begin{abstract}
This text is written in remembrance of Vladimir (Volodia) Sidorenko, honoring his scientific work, his kindness, and the warm and humorous spirit he brought into the community.
\end{abstract}

I was saddened and deeply shocked when I received the news from my colleague Violetta Weger that Volodia Sidorenko had passed away. I knew that his health had been fragile for some time and that he had been in the hospital on and off, but I never expected this outcome. The last email I received from him was in June—a warm message praising the memorial book I had compiled for another colleague, the late Ning Cai. It hurts to realize how final that message now is.

When I think of Volodia, I remember an exceptionally friendly, humorous, and helpful person. I hired him within two different research projects, and he always worked with great diligence and skill. Yet beyond his scientific contributions, it is his kindness, his small gestures, and his unique personality that stay with me most vividly.

As I try to write down my memories, I realize how little I truly know about the private corners of his life—and how much he still managed to share simply through his presence. Before joining my group, he had already worked at the ICE (Institute for Communications Engineering at Technical University of Munich). He quickly became known there as the colleague who would suddenly appear in the corridor or the kitchen with ice cream or chocolate for everyone. If I didn’t respond to his call to take a break, he would come to my office and insist—sometimes even placing the chocolate or ice cream directly on my desk. Small gestures, but remarkable in their warmth and consistency.

At conferences, he loved to walk around in the evenings offering vodka—``the best you can get in Russia,'' as he proudly emphasized. It was impossible not to enjoy those moments with him. You could talk to Volodia about science, of course, but also about everyday problems, family, worries, hopes. He had a calm, open way of listening and a gentle humor that made things feel lighter.

When my Russian friends Volodia Lebedev and Tanya Dolgova visited Munich, Volodia Sidorenko would always join us. He knew them well and enjoyed singing Russian songs with them late into the evening—something I remember with great fondness.

Scientifically, he entered my group through the NEWCOM project. Being an expert in coding theory, he began working on identification codes for me, in particular on the question of how useful random generators at the transmitter and receiver can be. Our collaboration had a very special rhythm: first I explained everything I knew—with many detailed questions from his side—and then I would hear nothing for quite a while… until he suddenly reappeared with an almost finished paper. We published that work \cite{Sidorenko2022} together, and I still smile when I recall that pattern.

After NEWCOM, he moved on to the QD-CAMNetz project and worked on quantum error correction. He took a great load off my shoulders by helping supervise one of my PhD students, Evagoras Stylianou, and he did it exceptionally well. Evagoras has told me many times how much he valued Volodia’s support, and I know he will miss him greatly—especially since my own departure from TUM at the end of 2023. Together with Volodia and Evagoras we published \cite{Stylianou2025,Stylianou2024,Sidorenko2023}.

Volodia was always honest, always open. I once played table tennis against him at a Christmas party—he kept insisting that his Chinese wife was much stronger than he was. I only met his child in person once, and then occasionally saw them appear in the background of Zoom meetings. His family will feel this loss most profoundly. He left them far too early.

Volodia’s passing leaves a painful emptiness—for his family, for his friends, and for all of us who had the joy of working with him. But I will remember him for his kindness, his humor, his modest brilliance, and all the small ways in which he made the world around him a warmer place.

\medskip

May this candle stand as a symbol of the light he brought into our lives.

He will not be forgotten.

\bigskip

\begin{flushright}
\textit{In deep gratitude and remembrance,}\\[1mm]
Christian Deppe
\end{flushright}

\end{document}